\documentclass[10pt]{amsart}

\usepackage{amssymb, amscd}
\usepackage[letterpaper,hypertex]{hyperref}
\usepackage{amsrefs}
\usepackage{diagrams}

\allowdisplaybreaks[1]

\theoremstyle{plain}
\newtheorem{theorem}{Theorem}
\newtheorem*{theorem*}{Theorem}

\newtheoremstyle{myexample*}{4pt}{4pt}{}{}{\bfseries}{.}{ }{\thmname{#1}:\thmnote{ #3}}

\theoremstyle{myexample*}
\newtheorem*{example*}{Example}

\theoremstyle{remark}
\newtheorem{remark}{Remark}

\newcommand{\RR}{\mathbb{R}}
\newcommand{\ZZ}{\mathbb{Z}}
\newcommand{\CC}{\mathbb{C}}

\newcommand{\PP}{\mathbb{P}}

\newcommand{\cl}{\operatorname{cl}}
\newcommand{\Cstar}{{\CC^\times}}

\newcommand{\ev}{\operatorname{ev}}
\newcommand{\Ev}{\operatorname{Ev}}

\newcommand{\Res}{\operatorname{Res}}
\newcommand{\lsp}{\mathfrak{sp}}

\DeclareMathOperator{\vdim}{vdim}

\newcommand{\cD}{\mathcal{D}}
\newcommand{\be}{\mathbf{e}}
\newcommand{\cF}{\mathcal{F}}

\newcommand{\cH}{\mathcal{H}}

\newcommand{\cL}{\mathcal{L}}

\newcommand{\bq}{\mathbf{q}}

\newcommand{\bt}{\mathbf{t}}

\renewcommand{\(}{\left(}
\renewcommand{\)}{\right)}

\begin{document}
\title[Givental's Lagrangian Cone]{Givental's Lagrangian Cone and
  $S^1$-Equivariant Gromov--Witten Theory}
\author{Tom Coates} 
\address{Department of Mathematics\\
  Harvard University\\
  1 Oxford Street\\
  Cambridge, MA 02138\\
  U.S.A.}
\curraddr{Department of Mathematics\\
  Imperial College London\\
  London, SW7 2AZ\\
  U.K.}
\email{tomc@imperial.ac.uk}
\thanks{This research was partially supported by the National Science
  Foundation grant DMS-0401275.}

\begin{abstract}
  In the approach to Gromov--Witten theory developed by Givental,
  genus-zero Gromov--Witten invariants of a manifold $X$ are encoded
  by a Lagrangian cone in a certain infinite-dimensional symplectic
  vector space.  We give a construction of this cone, in the spirit of
  $S^1$-equivariant Floer theory, in terms of $S^1$-equivariant
  Gromov--Witten theory of $X \times \PP^1$.  This gives a conceptual
  understanding of the ``dilaton shift'': a change-of-variables which
  plays an essential role in Givental's theory.
\end{abstract}

\maketitle
\section{Introduction}

It has long been understood that it is a good idea to arrange
Gromov--Witten invariants into generating functions which reflect
their origins in physics: many operations in Gromov--Witten theory
which seem complicated when viewed at the level of individual
invariants correspond to the application of rather simpler
differential operators to these generating functions.  A recent
insight of Givental is that such differential operators, which can
themselves appear quite complicated, are often the quantizations of
very simple linear symplectic transformations of a certain symplectic
vector space.  This point of view --- Givental's {\em quantization
  formalism} \cite{Givental:quantization, Givental:symplectic} --- has
been a crucial ingredient in several recent advances in the subject.
These include the proof of the Virasoro conjecture for toric Fano
manifolds \cite{Givental:semisimple}, the computation of twisted
Gromov--Witten invariants \cites{Coates--Givental:QRRLS, Tseng}, the
proof of a Quantum Hirzebruch--Riemann--Roch theorem relating quantum
extraordinary cohomology to quantum cohomology
\cites{Coates--Givental:Gelfand, Coates--Givental:QGC}, and the
construction of integrable hierarchies controlling the total
descendant potentials of certain Frobenius manifolds
\cites{Givental:nKdV, Givental--Milanov, Milanov}.

The symplectic vector space associated to the Gromov--Witten theory of
an almost-K\"ahler manifold $X$ is the space of Laurent series
\[
\cH = H^\bullet(X)\otimes \CC(\!(z^{-1})\!)
\]
equipped with the symplectic form
\[
\Omega(f,g) = \Res_{z=0} \Big(f(-z), g(z)\Big)\, dz.
\]
Here $(\cdot,\cdot)$ is the Poincar\'e pairing on $H^\bullet(X)$.
Generating functions for Gromov--Witten invariants --- the genus-$g$
Gromov--Witten potentials of $X$ and the total descendant potential of
$X$ --- are regarded as functions on $\cH_+ = H^\bullet(X)[z]$ via a
change of variables, called the {\em dilaton shift}, described in
equation \ref{eq:dilaton_shift} below.  Genus-zero Gromov--Witten
invariants are encoded by a certain Lagrangian submanifold $\cL$ of
$X$, defined in section \ref{section:genus_zero} below.  This submanifold $\cL$
has a very tightly-constrained geometry: it is a Lagrangian cone ruled
by a finite-dimensional family of isotropic subspaces
\cites{Coates--Givental:QRRLS,Givental: symplectic}.

We currently lack a conceptual understanding of why the quantization
formalism is so effective.  It makes sense, therefore, to look for a
geometric interpretation of the ingredients of the formalism --- of
the symplectic vector space $\cH$, the submanifold $\cL$, and the
dilaton shift.  In this paper we give a simple and geometric
construction of the submanifold $\cL$ in terms of the
$S^1$-equivariant Gromov--Witten theory of the space $X \times \PP^1$.
This gives rise to the dilaton shift in a natural way.  Our
construction suggests that $\cH$ should be thought of as the
$S^1$-equivariant Floer homology of the loop space of $X$; this is
discussed further in Section~\ref{section:localization} below.

The idea of the construction is as follows.  There is an ``evaluate at
infinity'' map
\[
\ev_\infty: (X \times \PP^1)^{op}_{0,n,(d,1)} \to X
\]
from an open set in the moduli space of stable maps to $X \times
\PP^1$ of bidegree $(d,1)$ from genus-zero curves with $n$ marked
points.  This open set is the locus of stable maps $f: \Sigma \to X
\times \PP^1$ such that the preimage $f^{-1}(X \times \{\infty\})$
consists of a single unmarked smooth point --- so there are no bubbles
or marked points over $\infty \in \PP^1$ --- and the map $\ev_\infty$
records the point of $X$ mapped to by $f^{-1}(X \times \{\infty\})$.
Although $\ev_\infty$ is not proper, it is equivariant with respect to
the $S^1$-action on $(X \times \PP^1)^{op}_{0,n,(d,1)}$ coming from
the $S^1$-action of weight $-1$ on the second factor of $X \times
\PP^1$ and the trivial $S^1$-action on $X$.  This allows us to define
a push-forward
\[
\(\ev_\infty\)_\star:
H^\bullet_{S^1}\((X \times \PP^1)^{op}_{0,n,(d,1)}\) \otimes \CC(\!(z^{-1})\!)
\to H^\bullet(X) \otimes \CC(\!(z^{-1})\!),
\]
where $H^\bullet_{S^1}(pt) = \CC[z]$: the restriction of the map
$\ev_\infty$ to $S^1$-fixed sets is proper, so we can define the
push-forward using fixed-point localization.  To push an equivariant
class forward along $\ev_\infty$ we first restrict it to the
$S^1$-fixed set in $(X \times \PP^1)^{op}_{0,n,(d,1)}$, then cap with
the virtual fundamental class of the fixed set, then divide by the
$S^1$-equivariant Euler class of the virtual normal bundle, and then
push forward (in the usual sense) along the map $\ev_\infty$ from the
$S^1$-fixed set to $X$.  One can think of this operation as a virtual
push-forward in $S^1$-equivariant Borel--Moore--Tate homology; it is
defined only over the field of fractions $\CC(z)$ of
$H^\bullet_{S^1}(pt)$, and not over $H^\bullet_{S^1}(pt)$ itself,
because we need to divide by the Euler class of the virtual normal
bundle.  The Lagrangian cone $\cL$ is the image of a certain class
\begin{equation}
  \label{eq:J-class}
  (-z) \sum_{\substack{d \in H_2(X;\ZZ) \\n \geq 0}} {Q^d \over n!}
  \prod_{i=1}^{i=n} \ev_i^\star \bt(\psi_i) \in
  \bigoplus_{\substack{d \in H_2(X;\ZZ) \\n \geq 0}}
  H^\bullet_{S^1}\((X \times \PP^1)^{op}_{0,n,(d,1)}\),
\end{equation}
defined in detail in Section \ref{section:quantization} below, under
this push-forward.

The dilaton shift arises here in the following way: the $S^1$-fixed
set in the space $(X \times \PP^1)^{op}_{0,n,(d,1)}$ can almost always
be identified with the space $X_{0,n+1,d}$ of degree-$d$ stable maps
to $X$ from genus-zero curves with $n+1$ marked points.  There are two
exceptions to this, however: when $(n,d) = (0,0)$ and when $(n,d) =
(1,0)$, the moduli space $X_{0,n+1,d}$ is empty but the $S^1$-fixed
set is a copy of $X$.  It is the contributions to the push-forward of
\eqref{eq:J-class} coming from these exceptional fixed loci which give
rise to the dilaton shift.  In the notation of Section
\ref{section:quantization}, the push-forward of \eqref{eq:J-class} is
\begin{multline*}
\underbrace{-z + \bt(z)}_{\substack{\text{contributions from} \\ \text{
      exceptional fixed loci}}} +\\
\sum_{\substack{d \in H_2(X;\ZZ) \\n \geq 0}} {Q^d \over n!}
\underbrace{\(\ev_{n+1}\)_\star \left[ [X_{0,n+1,d}]^{vir} \cap
  \(\prod_{i=1}^{i=n} \ev_i^\star \bt(\psi_i)\) \cdot
  {1 \over -z-\psi_{n+1}} \right]
}_{\text{contribution from $X_{0,n+1,d}$}}.
\end{multline*}
This makes the change of variables \eqref{eq:dilaton_shift} seem
very natural. 

We should emphasize that none of the geometric ingredients here are
new.  The observation that a product of two copies of the $J$-function
--- a certain generating function for genus-zero Gromov--Witten
invariants --- can be computed by fixed-point localization on the
graph space $X \times \PP^1$ was, or was equivalent to, a crucial step
in proving mirror theorems for toric varieties
\cite{Givental:equivariant, Lian--Liu--Yau:1, Lian--Liu--Yau:2,
  Lian--Liu--Yau:3, Givental:toric, Bertram}.  The equivariant
push-forward described above was introduced by Braverman
\cite{Braverman} in order to extract one copy of the $J$-function of a
flag manifold from the corresponding graph space.  The content of this
paper is the observation that when Braverman's construction is
extended to ``big quantum cohomology'' and to include gravitational
descendants, the dilaton shift emerges automatically.

Experts in the subject may wish to stop reading here, as what follows
is routine.  Section \ref{section:quantization} contains an
introduction to Givental's quantization formalism.  The details of the
construction of $\cL$ are in Theorem \ref{theorem:main} and Section
\ref{section:localization}.  The localization theorem which we need
does not appear to have been written down anywhere, so we prove it in
the Appendix.

\subsection*{Acknowledgements}
I would like to thank Alexander Braverman, who taught me the
construction on which this paper is based, and Mike Hopkins for
stimulating and useful discussions.  I am grateful also to the
Department of Mathematics at Imperial College London for hospitality
whilst this paper was being written.

\section{Givental's Quantization Formalism}
\label{section:quantization}

We begin by describing the quantization formalism.  We fix notation
for Gromov--Witten invariants in section \ref{section:GW} and discuss
the framework for working with higher-genus invariants in section
\ref{section:higher_genus}.  The latter section is not logically necessary:
the reader who is familiar with Givental's approach or uninterested in
the surrounding context should skip to section
\ref{section:genus_zero}, where the genus-zero picture is described.

\subsection{Gromov--Witten Invariants} \label{section:GW}

Let $X$ be a smooth projective variety\footnote{A virtual localization
  theorem has recently been established in the symplectic category
  \cite{Chen--Li, Hofer--Wysocki--Zehnder}, and so the constructions
  here can now be extended to the case of almost-K\"ahler target
  manifolds $X$.}.  The Gromov--Witten invariants of $X$ are certain
intersection numbers in moduli spaces of stable maps (see {\em e.g.}
\cites{Kontsevich--Manin, Ruan--Tian, Manin, Fulton--Pandharipande,
  Mirror}).  Let $X_{g,n,d}$ denote the moduli space of stable maps to
$X$ of degree $d \in H_2(X;\ZZ)$ from curves of genus $g$ with $n$
marked points, and let $[X_{g,n,d}]^{vir}$ be its virtual fundamental
class \cite{Li--Tian, Behrend--Fantechi, Behrend}.  The moduli space
comes equipped with evaluation maps
\begin{align*}
  \ev_i:X_{g,n,d} \to X && i\in\{1,\ldots,n\} \\
  \intertext{and universal cotangent line bundles}
  L_i \to X_{g,n,d} &&  i\in\{1,\ldots,n\}
\end{align*}
at each marked point.  We denote the first Chern class of $L_i$ by
$\psi_i$.  Gromov--Witten invariants are intersection numbers of the
form
\begin{equation}
  \label{eq:GW}
  \int_{[X_{g,n,d}]^{vir}} \prod_{i=1}^{i=n} \ev_i^\star(\alpha_i)
  \cdot \psi_i^{k_i},
\end{equation}
where $\alpha_1, \ldots, \alpha_n$ are cohomology classes on $X$ and
$k_1,\ldots,k_n$ are non-negative integers.  If any of the $k_i$ are
non-zero, such invariants are called gravitational descendants.

The {\em genus-$g$ descendant potential of $X$} is a generating
function for Gromov--Witten invariants:
\[
\cF^g_X(t_0,t_1,\ldots) =
\sum_{d \in H_2(X;\ZZ)} \sum_{n \geq 0} {Q^d \over n!}
\int_{[X_{g,n,d}]^{vir}} \prod_{i=1}^{i=n} \ev_i^\star \bt(\psi_i).
\]
Here $t_0, t_1,\ldots$ are cohomology classes on $X$; $\bt(\psi) = t_0
+ t_1 \psi + t_2 \psi^2 + \ldots$, so that
\begin{equation}
  \label{eq:vector_pullback}
  \ev_i^\star \bt(\psi_i) = \ev_i^\star(t_0) +
  \ev_i^\star(t_1) \cdot \psi_i + \ev_i^\star (t_2) \cdot \psi_i^2 +
  \ldots;
\end{equation}
and $Q^d$ is the representative of $d$ in the Novikov ring \cite[III
5.2.1]{Manin}, which is a certain completion of the group ring of
$H_2(X;\ZZ)$.  If we pick a basis $\{\phi_1,\ldots,\phi_N\}$ for
$H^\bullet(X;\CC)$ and write
\begin{equation}
  \label{eq:co-ords}
  t_i = t_i^1 \phi_1 + \ldots + t_i^N \phi_N
\end{equation}
then
\[
\cF^g_X(t_0,t_1,\ldots) =
\sum_{\substack{d \in H_2(X;\ZZ) \\n \geq 0}}
\sum_{\substack{k_1,\ldots,k_n \\ \alpha_1,\ldots,\alpha_n}}
{Q^d t_{k_1}^{\alpha_1} \ldots t_{k_n}^{\alpha_n} \over n!}
\int_{[X_{g,n,d}]^{vir}} \prod_{i=1}^{i=n}
\ev_i^\star(\phi_{\alpha_i}) \cdot \psi_i^{k_i},
\]
so we can regard $\cF^g_X$ as a formal power series with Taylor
coefficients given by Gromov--Witten invariants of $X$.  The {\em
  total descendant potential of $X$}
\[
\cD_X(t_0,t_1,\ldots) = \exp \( \sum_{g \geq 0} \hbar^{g-1}
\cF^g_X(t_0,t_1,\ldots) \)
\]
is a generating function for Gromov--Witten invariants of all genera.

\subsection{The quantization formalism} \label{section:higher_genus}

Consider the space
\[
\cH = H^\bullet(X)(\!(z^{-1})\!)
\]
of cohomology-valued Laurent series, equipped with the symplectic form
\[
\Omega(f,g) = \Res_{z=0} \Big( f(-z), g(z) \Big) \, dz.
\]
Here and from now on we work over a ground ring $\Lambda$ which is the
tensor product of the Novikov ring with $\CC$: we take cohomology with
coefficients in $\Lambda$, the Poincar\'e pairing $(\cdot,\cdot)$ and
the symplectic form are $\Lambda$-valued, and so on.  The space $\cH$
is the direct sum of Lagrangian subspaces
\begin{align*}
  \cH_+ = H^\bullet(X)[z], && \cH_- =
  z^{-1} H^\bullet(X)[\![z^{-1}]\!].
\end{align*}
A general element of $\cH$ takes the form
\begin{equation}
  \label{eq:Darboux}
  \sum_{k=0}^{\infty} \sum_{\mu = 1}^{N} q^\mu_k \, \phi_\mu z^k +
  \sum_{l=0}^\infty \sum_{\nu=1}^N p^\nu_l \, \phi^\nu (-z)^{-1-l}
\end{equation}
where $\{\phi_1,\ldots,\phi_N\}$ is our basis for $H^\bullet(X)$, we
set $g_{\alpha \beta} = \(\phi_\alpha,\phi_\beta\)$, define $g^{\alpha
  \beta}$ to be the $(\alpha,\beta)$-entry of the matrix inverse to
that with $(\alpha,\beta)$-entry $g_{\alpha \beta}$, and raise indices
using $g^{\alpha \beta}$:
\[
\phi^\nu = \sum_{\lambda = 1}^N g^{\nu \lambda} \phi_\lambda.
\]
Equation \eqref{eq:Darboux} defines Darboux co-ordinates
$\{q^\mu_k, p^\nu_l\}$ on $\cH$ which are compatible with the
polarization $\cH = \cH_+ \oplus \cH_-$.

To each linear infinitesimal symplectic transformation $A \in
\lsp\(\cH\)$ we associate a differential operator --- the {\em
  quantization} of $A$ --- constructed as follows.  The quadratic
Hamiltonian $h_A:f \mapsto {1 \over 2} \Omega(Af,f)$ can be written as
a linear combination of quadratic monomials in the Darboux
co-ordinates $\{q^\mu_k, p^\nu_l\}$.  We set
\begin{align*}
  \widehat{q^\mu_k q^\nu_l} = {q^\mu_k q^\nu_l \over \hbar}, &&
  \widehat{p^\mu_k q^\nu_l} = q^\nu_l {\partial \over \partial
  q^\mu_k}, &&
  \widehat{p^\mu_k p^\nu_l} = \hbar {\partial \over \partial q^\mu_k} {\partial \over \partial q^\nu_l},
\end{align*}
and extend by linearity, defining the quantization $\widehat{A}$ of
$A$ to equal $\widehat{h_A}$.  The quantized operator $\widehat{A}$
acts on certain\footnote{Since the symplectic space $\cH$ is
  infinite-dimensional, quantizations $\widehat{A}$ in general contain
  infinite sums of differential operators.  The application of a
  general quantized infinitesimal symplectic transformation to a
  general formal power series in the variables $q^\alpha_k$ is not
  well-defined, but it is easy to check that the operations used in
  the Example below do in fact make sense.} formal power series in the
variables $q^\alpha_k$, where $\alpha \in \{ 1,\ldots, N\}$ and $k
\geq 0$.

Let
\begin{align*}
  q_k &= \sum_{\lambda = 1}^N q^\lambda_k \phi_\lambda &
  k=0,1,2,\ldots,
\end{align*}
and
\[
\bq(z) = q_0 + q_1 z + q_2 z^2 + \ldots.
\]
Quantized infinitesimal symplectic transformations $\widehat{A}$ act
on certain formal functions of $\bq(z)$ --- in other words, on
certain formal power series in the $q^\alpha_k$ --- whereas the
total descendant potential $\cD_X(t_0,t_1,\ldots)$ is a formal
function of
\[
\bt(z) = t_0 + t_1 z + t_2 z^2 + \ldots
\]
--- or in other words, a formal power series in the variables
$t^\alpha_k$ from \eqref{eq:co-ords}.  We let quantized operators
$\widehat{A}$ act on the total descendant potential
$\cD_X(t_0,t_1,\ldots)$ via the identification
\begin{equation}
  \label{eq:dilaton_shift}
  \bq(z) = \bt(z) - z.
\end{equation}
This change of variables is called the {\em dilaton shift}.

This framework allows one to express many operations which arise in
Gromov--Witten theory in terms of the quantizations of very simple
linear symplectic transformations of $\cH$.  One example of this
occurs in the Gromov--Witten theory of a point.

\begin{example*}[The Virasoro Conjecture]
  Let $X$ be a point.  The corresponding symplectic space is
  \begin{align*}
    \cH = \CC(\!(z^{-1})\!), && \Omega(f,g) = \Res_{z=0} f(-z) g(z) \, dz,
  \end{align*}
  and Darboux co-ordinates $\{q_k, p_l\}$ on $(\cH,\Omega)$ are given
  by
  \[
  \ldots {p_2 \over (-z)^3} + {p_1 \over (-z)^2} + {p_0 \over (-z)} +
  q_0 + q_1 z + q_2 z^2 + \ldots.
  \]
  The quadratic Hamiltonians corresponding to the linear infinitesimal
  symplectic transformations
\begin{align*}
  l_n: \cH & \longrightarrow \cH \\
  f & \longmapsto z^{-1/2} \(z {d \over d z} z\)^{n+1} z^{-1/2} f && n
  \geq -1
\end{align*}
are
\begin{align*}
  & - \sum_{k \geq 1} p_{k-1} q_k -{1 \over 2} q_0^2 & n=-1 \\
  & - \sum_{k \geq 0} {\Gamma(k+n+{3 \over 2}) \over \Gamma(k+{1 \over
      2})} q_k p_{k+n} + \sum_{l=0}^{l=n-1} (-1)^l {\Gamma(n-l+{1
      \over 2}) \over \Gamma(-l-{1 \over 2})} p_l p_{n-1-l} & n \geq
  0
\end{align*}
The quantizations $\widehat{l_n}$ are the differential operators
\begin{align*}
   {\partial \over \partial t_0} & - \sum_{k \geq 1} t_k
  {\partial \over \partial t_{k-1}} + {t_0^2 \over 2 \hbar} & n = -1 \\
  {\Gamma\(n + {5 \over 2}\) \over \Gamma\({3 \over
      2}\)} {\partial \over \partial t_{n+1}} & -
  \sum_{k \geq 0} {\Gamma\(k+n+{3 \over 2}\) \over \Gamma\(k + {1
      \over 2}\)}  t_k {\partial \over \partial t_{k+n}} \\
  & - {\hbar \over 2}
  \sum_{l=0}^{l=n-1} (-1)^{l+1} {\Gamma\(n-l+{1 \over 2}\) \over \Gamma\(-l- {1
      \over 2}\)} {\partial \over \partial t_l} {\partial \over
    \partial t_{n-1-l}}
  & n \geq 0.
\end{align*}
Note that the dilaton shift \eqref{eq:dilaton_shift} plays an
essential role here, as without it these differential operators would
not be quadratic in $p_\alpha$ and $q_\beta$.   The Virasoro Conjecture
for Gromov--Witten invariants of a point (see {\em e.g.}
\cite{Getzler:Virasoro}) asserts that
\begin{align*}
\( \widehat{l_n} - {\delta_{n,0} \over 16}\)
\cD_{pt}(t_0,t_1,\ldots) = 0, && n \geq -1.
\end{align*}
This is equivalent \cite{Dijkgraaf--Verlinde--Verlinde:Virasoro} to
Witten's Conjecture \cite{Witten:KdV}, proved by Kontsevich
\cite{Kontsevich:KdV}.
\end{example*}

\subsection{The genus-zero picture} \label{section:genus_zero}

So far we have considered a formalism for working with Gromov--Witten
invariants of all genera.  This involves quantized symplectic
transformations applied to generating functions for the invariants.
The semi-classical limit of this framework involves {\em unquantized}
symplectic transformations applied to certain Lagrangian submanifolds
of $\cH$.  This is how the Lagrangian submanifold $\cL$ from the
Introduction enters the theory.

It is easy to see that if
\[
\cD(s) = \exp\(\sum_{g \geq 0} \hbar^{g-1} \cF^g(s)\)
\]
is a one-parameter family of formal power series in the variables
$q^\mu_k$ such that
\[
{d \over ds} \cD(s) = \widehat{A} \, \cD(s)
\]
for some $A \in \lsp\(\cH\)$, then the formal germ of a Lagrangian
submanifold of $\cH$ given in Darboux co-ordinates \eqref{eq:Darboux}
by
\[
p^\nu_l = {\partial \cF^0(s) \over \partial q^\nu_l}
\]
evolves with $s$ under the Hamiltonian flow of $h_A$.  We thus
consider the formal germ of a Lagrangian submanifold $\cL$ defined by
\begin{equation}
  \label{eq:cone}
  p^\nu_l = {\partial \cF_X^0 \over \partial q^\nu_l},
\end{equation}
where we regard $\cF^0_X(t_0,t_1,\ldots)$ as a formal power series in
the $q^\nu_l$ via the dilaton shift \eqref{eq:dilaton_shift}.  The
formal germ $\cL$ is defined for $\bq(z)$ near $-z$.  It corresponds,
under the identification of $\cH = \cH_+ \oplus \cH_-$ with $T^\star
\cH_+ = \cH_+ \oplus \cH_+^\vee$ coming from the polarization, to the
graph of the differential of the genus-zero descendant potential
$\cF^0_X$.  $\cL$ therefore encodes genus-zero Gromov--Witten
invariants of $X$.  A general point of $\cL$ takes the form
\begin{multline}
  \label{eq:better_cone}
  \bq(z) \, + \\\sum_{\substack{d \in H_2(X;\ZZ) \\n \geq 0}} {Q^d \over
    n!} \(\ev_{n+1}\)_\star \left[ [X_{0,n+1,d}]^{vir} \cap
    \(\prod_{i=1}^{i=n} \ev_i^\star
    \bt(\psi_i)\) \cdot {1 \over - z - \psi_{n+1}} \right].
\end{multline}
To see this, expand ${1 \over - z - \psi_{n+1}}$ as a power series in
$z^{-1}$ and compare \eqref{eq:better_cone} with \eqref{eq:Darboux} and
\eqref{eq:cone}.

\section{The Localization Calculation}
\label{section:localization}

We begin this section by giving a precise definition of the virtual
push-forward described in the Introduction.  We then state Theorem 1.
The proof of Theorem 1, which is a straightforward application of the
virtual localization result of Graber and Pandharipande
\cite{Graber--Pandharipande}, is contained in section
\ref{section:proof}.

\subsection{A virtual push-forward} \label{section:push-forward}

Given schemes $Y$ and $Z$ with $\Cstar$-action\footnote{We have
  switched from $S^1$-actions to $\Cstar$-actions in order to make use
  of the virtual localization result \cite{Graber--Pandharipande}.},
an equivariant map $f:Y \to Z$ such that the induced map on fixed sets
is proper gives a push-forward
\begin{equation}
  \label{eq:push}
  f_\star: H^\Cstar_{\bullet,BM}(Y)\otimes \CC(z) \to
  H^\Cstar_{\bullet,BM}(Z)\otimes \CC(z)
\end{equation}
in $\Cstar$-equivariant\footnote{Equivariant Borel--Moore homology is
  discussed in the Appendix.} Borel--Moore homology \cite{Braverman}.
$\CC(z)$ here is the field of fractions of $H^\bullet_\Cstar(pt) =
\CC[z]$.  The localization theorem (see
\cite{Edidin--Graham:localization} and the Appendix) implies that the
maps
\begin{align*}
  \(i_Y\)_\star:H^\Cstar_{\bullet,BM}(Y^\Cstar) \to
  H^\Cstar_{\bullet,BM}(Y), &&
  \(i_Z\)_\star:H^\Cstar_{\bullet,BM}(Z^\Cstar) \to
  H^\Cstar_{\bullet,BM}(Z)
\end{align*}
induced by the inclusions $i_Y:Y^\Cstar \to Y$, $i_Z:Z^\Cstar \to Z$
of $\Cstar$-fixed sets become isomorphisms after tensoring with
$\CC(z)$.  The push-forward \eqref{eq:push} is defined to be the
composition
\begin{diagram}
  H^\Cstar_{\bullet,BM}(Y) \otimes \CC(z) &
  \rDashto^{f_\star} &
  H^\Cstar_{\bullet,BM}(Z) \otimes \CC(z) \\
  \dTo^{\(\(i_Y\)_\star\)^{-1}} & &
  \uTo_{\(i_Z\)_\star} \\
  H^\Cstar_{\bullet,BM}(Y^\Cstar) \otimes \CC(z) &
  \rTo &
  H^\Cstar_{\bullet,BM}(Z^\Cstar) \otimes \CC(z)
\end{diagram}
where the bottom horizontal arrow is the usual proper push-forward.
When the map $f$ is proper, \eqref{eq:push} agrees with the usual
push-forward.

If $Y$ and $Z$ are smooth $\Cstar$-varieties and $f:Y \to Z$ is
equivariant and proper on fixed sets as before then this construction
gives a push-forward in equivariant cohomology
\[
f_\star:H^\bullet_\Cstar(Y) \otimes \CC(z) \to H^\bullet_\Cstar(Z) \otimes
\CC(z)
\]
which raises degree by $2 \dim_\CC(Z) - 2\dim_\CC(Y)$.  This is by definition
the composition
\begin{diagram}
  H^\bullet_\Cstar(Y) \otimes \CC(z) & \rDashto^{f_\star} &
  H^{\bullet}(Z) \otimes \CC(z) \\
  \dTo & & \uTo \\
  H^\Cstar_{\bullet, BM}(Y) \otimes \CC(z) & \rTo &
  H^\Cstar_{\bullet,BM}(Z) \otimes \CC(z)
\end{diagram}
where the vertical arrows are Poincar\'e duality and the bottom
horizontal arrow is the push-forward \eqref{eq:push}.

In the case we wish to consider, $Y$ will be an open subset of a
moduli space of stable maps.  This need not be smooth, but it it does
carry a $\Cstar$-equivariant perfect obstruction theory: it is
``virtually smooth''.  Given a $\Cstar$-scheme $Y$ equipped with a
$\Cstar$-equivariant perfect obstruction theory, a smooth
$\Cstar$-variety $Z$, and an equivariant map $f:Y \to Z$ which is
proper on fixed sets, we define the {\em virtual push-forward}
\[
f_\star:H^\bullet_\Cstar(Y) \otimes \CC(z) \to H^\bullet_\Cstar(Z) \otimes
\CC(z)
\]
as follows.  The obstruction theory determines a virtual fundamental
class \cite{Li--Tian, Behrend--Fantechi, Behrend} in the equivariant
Chow group $A^\Cstar_{\vdim(Y)}(Y)$, where $\vdim(Y)$ is the virtual
dimension, and hence (via the cycle map) gives a class in equivariant
Borel--Moore homology
\[
[Y]^{vir} \in H_{2 \vdim(Y), BM}^\Cstar(Y).
\]
The virtual push-forward is defined to be the composition
\begin{diagram}
  H^\bullet_\Cstar(Y) \otimes \CC(z) & \rDashto^{f_\star} &
  H^{\bullet}(Z) \otimes \CC(z) \\
  \dTo^{[Y]^{vir} \cap \ } & & \uTo \\
  H^\Cstar_{\bullet, BM}(Y) \otimes \CC(z) & \rTo &
  H^\Cstar_{\bullet,BM}(Z) \otimes \CC(z)
\end{diagram}
where the left-hand vertical arrow is cap product with the class
$[Y]^{vir}$, the right-hand vertical arrow is Poincar\'e duality, and
the bottom horizontal arrow is the push-forward \eqref{eq:push}.  The
virtual push-forward raises degree by $2 \dim_\CC(Z) - 2\vdim(Y)$.
Once appropriate definitions are in place, the construction extends
without change to the case (which we will need below) in which $Y$ is
a Deligne--Mumford quotient stack rather than a scheme; see the
Appendix for details.

The virtual localization result of Graber and Pandharipande
\cite{Graber--Pandharipande} implies that, under a mild technical
hypothesis\footnote{They require that $Y$ admit a $\Cstar$-equivariant
  embedding into a non-singular Deligne--Mumford stack.  This is the
  case when $Y$ is an open subset of a moduli stack of stable maps to
  a $\Cstar$-variety: see Appendices A and C of
  \cite{Graber--Pandharipande}.},
\begin{equation}
  \label{eq:virtual_localization}
  [Y]^{vir} = \(i_Y\)_\star \left[ \sum {[Y_{j}]^{vir} \over
  \be\(N_{j}^{vir}\)} \right] \in H^\Cstar_{\bullet,BM}(Y) \otimes \CC(z).
\end{equation}
The sum here is over components $Y_{j}$ of the $\Cstar$-fixed locus in
$Y$.  The virtual fundamental classes $[Y_{j}]^{vir}$ and virtual
normal bundles $N_{j}^{vir}$ are determined by the obstruction theory;
$\be$ here denotes the $\Cstar$-equivariant Euler class.  If we write
$f_{j}$ for the restriction of $f:Y \to Z$ to the $\Cstar$-fixed
component $Y_{j}$ then \eqref{eq:virtual_localization} implies that we
can write the virtual push-forward of $\alpha \in
H^\bullet_\Cstar(Y)\otimes \CC(z)$ as
\begin{equation}
  \label{eq:formula}
  f_\star(\alpha) =
  \sum \(f_{j}\)_\star \left[ \,
    {[Y_{j}]^{vir} \cap \left. \alpha \right|_{Y_{j}} \over
      \be(N_{j}^{vir})} \, \right].
\end{equation}

Consider now the open subset $(X \times \PP^1)^{op}_{0,n,(d,1)}$ of
the moduli space $(X \times \PP^1)_{0,n,(d,1)}$ consisting of those
stable maps $f: \Sigma \to X \times \PP^1$ such that the preimage
$f^{-1}(X \times \{\infty\})$ is a single unmarked smooth point
$x_\infty$.  Consider the $\Cstar$-action on moduli space coming from
the trivial $\Cstar$-action on $X$ and the $\Cstar$-action of weight
$-1$ on $\PP^1$.  The space $(X \times \PP^1)_{0,n,(d,1)}$ carries a
canonical $\Cstar$-equivariant perfect obstruction theory, and so the
$\Cstar$-invariant open subset $(X \times \PP^1)^{op}_{0,n,(d,1)}$
does too.  The ``evaluate at infinity'' map
\[
\ev_\infty: (X \times \PP^1)^{op}_{0,n,(d,1)}  \longrightarrow X
\]
which sends the stable map $f: \Sigma \to X \times \PP^1$ to
$f(x_\infty)$ is $\Cstar$-equivariant and proper on fixed sets.  The
virtual push-forwards along the maps $\ev_\infty$ assemble to give a
map
\[
\Ev_\infty: \bigoplus_{\substack{d \in H_2(X;\ZZ) \\n \geq 0}}
H^\bullet_{S^1}\((X \times \PP^1)^{op}_{0,n,(d,1)}\)
\longrightarrow \cH.
\]
We are now ready to state our result.

\begin{theorem} \label{theorem:main}
  $\cL$ is the image under $\Ev_\infty$ of the class
  \[
  (-z) \sum_{\substack{d \in H_2(X;\ZZ) \\n \geq 0}} {Q^d \over n!}
  \prod_{i=1}^{i=n} \ev_i^\star \bt(\psi_i) \in
  \bigoplus_{\substack{d \in H_2(X;\ZZ) \\n \geq 0}}
  H^\bullet_{S^1}\((X \times \PP^1)^{op}_{0,n,(d,1)}\).
  \]
\end{theorem}

\subsection{The Proof of Theorem 1} \label{section:proof}

This is a straightfoward application of the formula \eqref{eq:formula}
for the virtual push-forward.  The calculations are similar to, but
easier than, those occurring in section 4 of
\cite{Graber--Pandharipande}.

\subsubsection*{Case 1: $(n,d) \not \in \{(0,0), (1,0)\}$}  The $\Cstar$-fixed
locus in $(X \times \PP^1)^{op}_{0,n,(d,1)}$ consists of stable maps
from nodal curves such that exactly one component of the curve is
mapped with degree $1$ to $\{x_\infty\} \times \PP^1 \subset X \times
\PP^1$, and the rest of the curve is mapped to $X \times \{0\}$.  We
identify the fixed locus with the moduli space $X_{0,n+1,d}$ of
$(n+1)$-pointed stable maps to $X$: the component mapped to
$\{x_\infty\} \times \PP^1$ is attached at the $(n+1)$st marked point.
The $\Cstar$-fixed part of the perfect obstruction theory on $(X
\times \PP^1)^{op}_{0,n,(d,1)}$ coincides with the usual perfect
obstruction theory on $X_{0,n+1,d}$, and the virtual normal bundle to
the fixed locus is
\[
\CC_{(-1)} \oplus \(L_{n+1} \otimes \CC_{(-1)}\)
\]
where $\CC_{(-1)}$ denotes the trivial bundle over $X_{0,n+1,d}$ with
$\Cstar$-weight $-1$.  Thus
\begin{multline} \label{eq:general}
  \(\ev_\infty\)_\star \left[
    (-z) \prod_{i=1}^{i=n} \ev_i^\star \bt(\psi_i) \right] = \\
  \(\ev_{n+1}\)_\star \left[
    [X_{0,n+1,d}]^{vir} \cap
    \(\prod_{i=1}^{i=n} \ev_i^\star \bt(\psi_i)\) \cdot
    {1 \over - z - \psi_{n+1}} \right]
\end{multline}

\subsubsection*{Case 2: $(n,d) = (1,0)$}  We have
\[
\(X \times \PP^1\)^{op}_{0,n,(d,1)} \cong X \times \CC
\]
and the $\Cstar$-fixed locus here is a copy of $X$.  The virtual
fundamental class on $X$ determined by the $\Cstar$-fixed part of the
perfect obstruction theory is the usual fundamental class of $X$.  The
restriction to the fixed locus of the dual to the universal cotangent
line bundle $L_1$ is the trivial bundle $\CC_{(-1)}$ over $X$ of
$\Cstar$-weight $-1$, and the virtual normal bundle is also
$\CC_{(-1)}$.  Thus
\begin{equation}
  \label{eq:1,0}
  \(\ev_\infty\)_\star \Big[
  (-z) \cdot \ev_1^\star \bt(\psi_1) \Big] = \bt(z).
\end{equation}

\subsubsection*{Case 3: $(n,d) = (0,0)$}  Here
\[
\(X \times \PP^1\)^{op}_{0,0,(0,1)} \cong X
\]
and there is no moving part of the obstruction theory.  The virtual
fundamental class induced on the fixed locus $X$ is the usual
fundamental class of $X$, and
\begin{equation}
  \label{eq:0,0}
  \(\ev_\infty\)_\star \Big[ {-z} \Big] = -z.
\end{equation}

\medskip

Combining \eqref{eq:general}, \eqref{eq:1,0}, and
\eqref{eq:0,0}, we find that the image of the class from Theorem 1
under $\Ev_\infty$ is
\begin{multline*}
  -z + \bt(z) + \\\sum_{\substack{d \in H_2(X;\ZZ) \\n \geq 0}} {Q^d \over
    n!} \(\ev_{n+1}\)_\star \left[ [X_{0,n+1,d}]^{vir} \cap
    \(\prod_{i=1}^{i=n} \ev_i^\star
    \bt(\psi_i)\) \cdot {1 \over - z - \psi_{n+1}} \right].
\end{multline*}
This coincides with our expression \eqref{eq:better_cone} for a
general point of $\cL$.  The proof is complete.  \hfill $\qed$

\begin{remark}
  We see from the proof of Theorem \ref{theorem:main} that one should
  regard the factor of $-z$ occurring in the statement as the
  $\Cstar$-equivariant Euler class of $R^\bullet \pi_\star
  \ev_{n+1}^\star \CC_{(-1)}$, where $\pi:X_{g,n+1,d} \to X_{g,n,d}$
  is the universal family over the moduli space of stable maps and
  $\CC_{(-1)}$ is the trivial bundle of $\Cstar$-weight $-1$ over $X$.
  Such a ``twist by the Euler class'' roughly corresponds to
  considering the Gromov--Witten theory of a hypersurface
  \cite{Coates--Givental:QRRLS}.  If we regard our study of $\(X
  \times \PP^1\)^{op}_{0,n,(d,1)}$ as a proxy for studying the
  Gromov--Witten theory of $X \times \CC$ then the two ingredients of
  our construction push in opposite directions: we end up, roughly
  speaking, thinking of $X$ as an ``equivariant hypersurface'' in $X
  \times \CC_{(-1)}$.  The dilaton shift arises exactly from the
  difference between the two notions of stability here: stability as a
  map to $X$ and stability as a graph in $X \times \CC$.
\end{remark}

\begin{remark} \label{Floerremark} Our construction of $\cL$ bears a
  striking resemblance to the ``fundamental Floer cycle'' --- the
  semi-infinite cycle in loop space consisting of loops which bound
  holomorphic discs --- in the heuristic picture relating quantum
  cohomology to the $S^1$-equivariant Floer homology of loop space
  outlined in \cite{Givental:heuristic}.  This suggests that one
  should regard $\cH$ as the $S^1$-equivariant Floer homology of the
  loop space of $X$.  Other evidence for this comes from comparing the
  symplectic transformation in \cite[Theorem
  1]{Coates--Givental:QRRLS} with the calculations in \cite[Section
  4]{Givental:heuristic}, and from the beautiful recent work of
  Costello \cite{Costello}.  As mentioned above, the graph space $\(X
  \times \PP^1\)_{0,n,(d,1)}$ plays a key role in many proofs of toric
  mirror symmetry \cite{Givental:equivariant,
    Lian--Liu--Yau:1,Lian--Liu--Yau:2,Lian--Liu--Yau:3,
    Givental:toric, Kim, Bertram, Lee}, where it links Floer-theoretic
  predictions to rigorous calculations in Gromov--Witten theory.  It
  would be interesting to understand exactly how $S^1$-equivariant
  Floer homology relates to our picture.
\end{remark}

\section*{Appendix: $\Cstar$-Equivariant Borel--Moore Homology}

In \cite{Braverman} Braverman used a sheaf-theoretic definition of
equivariant Borel--Moore homology, in the spirit of \cite{Lusztig}. We
will take a different point of view, regarding Borel--Moore homology
as the homology theory of singular chains with locally finite support.
This meshes more readily with constructions of the virtual fundamental
class.  We collect the properties of non-equivariant Borel--Moore
homology that we will need in section A1 and describe the equivariant
theory, constructed by Edidin and Graham in
\cite{Edidin--Graham:intersection}, in section A2.  In section A3 we
discuss the Borel--Moore homology of certain quotient stacks.  Since
the precise form of the localization theorem for $\Cstar$-equivariant
Borel--Moore homology which we used in section
\ref{section:push-forward} does not appear to have been written down
anywhere, we prove it in section A4; it was undoubtedly already
well-known.

\subsection*{A1. Borel--Moore homology}

Good introductions to Borel--Moore homology can be found in
\cite[chapter 19]{Fulton:intersection}, \cite[Appendix
B]{Fulton:Young}, \cite[section 2.6]{Chriss--Ginzburg}, and
\cite[Appendix C]{Nakajima--Yoshioka}.  We work with the definition
from \cite{Fulton:intersection}: if a space $X$ is embedded as a
closed subspace of $\RR^n$ then
\begin{equation}
  \label{eq:BM}
  H_{i,BM}(X) := H^{n-i}(\RR^n, \RR^n-X).
\end{equation}
All homology and cohomology groups are taken with complex coefficients
throughout.  Properties of Borel--Moore homology include:
\begin{itemize}
\item[{\bf BM1}] There are {\em cap products}
  \[
  H^j(X) \otimes H_{k,BM}(X) \to H_{k-j,BM}(X).
  \]
  See \cite[section 19.1]{Fulton:intersection}.
\item[{\bf BM2}] If $X$ is a smooth variety of dimension $n$ then
  $H_{2n,BM}(X)$ is freely generated by the {\em fundamental class}
  $[X] \in H_{2n,BM}(X)$, and
  \[
  [X] \cap: H^k(X) \to H_{2n-k,BM}(X)
  \]
  is an isomorphism.  This is {\em Poincar\'e duality}.  See
  \cite[section 19.1]{Fulton:intersection}.
\item[{\bf BM3}] There is a {\em K\"unneth formula}
  \[
  H_{k,BM}(X \times Y) = \bigoplus_{i+j=k} H_{i,BM}(X) \otimes
  H_{j,BM}(Y).
  \]
  This follows immediately from definition \eqref{eq:BM} and the
  K\"unneth formula for relative homology.
\item[{\bf BM4}] There are {\em covariant push-forwards} for proper
  maps $f:X \to Y$,
  \[
  f_\star:H_{k,BM}(X) \to H_{k,BM}(Y).
  \]
  See   \cite[section 19.1]{Fulton:intersection}.
\item[{\bf BM5}] There are {\em contravariant pull-backs} for open
  embeddings $j:U \to Y$,
  \[
  j^\star:H_{k,BM}(Y) \to H_{k,BM}(U).
  \]
  See   \cite[section 19.1]{Fulton:intersection}.
\item[{\bf BM6}] There is a {\em long exact sequence}
  \[
  \ldots \to H_{i+1,BM}(U) \to H_{i,BM}(X) \xrightarrow{i_\star}
  H_{i,BM}(Y) \xrightarrow{j^\star} H_{i,BM}(U) \to \ldots
  \]
  where $j:U \to Y$ is an open embedding and $i:X \to Y$ is the closed
  embedding of the complement $X$ to $U$ in $Y$.    See
  \cite[section 19.1]{Fulton:intersection}.
\item[{\bf BM7}] If $X$ is a scheme of dimension $n$ then
  $H_{i,BM}(X)=0$ for $i>2n$.  This is part of Lemma 19.1.1 in
  \cite{Fulton:intersection}.
\item[{\bf BM8}] For any scheme $X$ there is a {\em cycle map}
  \[
  \cl:A_k(X) \to H_{2k,BM}(X)
  \]
  which is covariant for proper maps and compatible with Chern
  classes.  See \cite[section 19.1]{Fulton:intersection}.
\item[{\bf BM9}] For any l.c.i. morphism of schemes $f:Y \to X$ of
  relative dimension $d$ there is a {\em Gysin map}
  \[
  f^\star:H_{k,BM}(X) \to H_{k-2d, BM}(Y).
  \]
  Such maps are functorial and compatible with the cycle class.  When
  $Y$ is a vector bundle over $X$ of rank $d$, $f^\star$ is the {\em
    Thom isomorphism} $H_{k,BM}(X) \to H_{k+2d,BM}(Y)$.  See
  \cite[Example 19.2.1]{Fulton:intersection}.

\end{itemize}

\subsection*{A2. Equivariant Borel--Moore homology}

Given a $g$-dimensional linear algebraic group $G$ acting in a
reasonable way\footnote{We sidestep a technical issue here.  Edidin
  and Graham work with algebraic spaces, rather than schemes.  This is
  because the quotient of an algebraic space by a free action of an
  algebraic group is an algebraic space, but the quotient of a scheme
  by a free action of of an algebraic group need not be a scheme.  We
  would like the mixed space $X_G$ to be a scheme, because we want to
  use properties of the Borel--Moore homology of schemes listed in
  section A1.  Proposition 23 in \cite{Edidin--Graham:intersection}
  gives conditions on the group action sufficient to ensure that $X_G$
  is a scheme: we will consider only actions of $G$ on $X$ which
  satisfy these hypotheses, calling such actions {\em reasonable}.  In
  view of the construction of the moduli space of stable maps as a
  stack quotient given in \cite{Fulton--Pandharipande}, it suffices
  for the purposes of this paper to consider only reasonable actions.
  Another, perhaps more satisfactory, approach would be to develop a
  Borel--Moore homology theory for algebraic spaces --- much as is
  done for intersection theory in section 6.1 of
  \cite{Edidin--Graham:intersection} --- but as we do not need to do
  this, we won't.} on an scheme $X$ of dimension $n$, Edidin and
Graham \cite{Edidin--Graham:intersection} define the $G$-equivariant
Borel--Moore homology groups of $X$ as
\[
H^G_{i,BM}(X) := H_{i+2l-2g,BM}(X_G).
\]
Here $X_G$ is the {\em mixed space} $(X \times U)/G$, where $U$ is an
open set in an $l$-dimensional representation $V$ of $G$ such that the
action of $G$ on $U$ is free and the real codimension of $V-U$ in $V$
is more than $2n-i+1$.

One can see that this is well-defined using Bogomolov's {\em double
  filtration argument} \cite[Definition-Proposition 1 and Section
2.8]{Edidin--Graham:intersection}.  Suppose that $V_1$ and $V_2$ are
representations of $G$ respectively of dimensions $l_1$ and $l_2$ and
containing open sets $U_1$ and $U_2$ such that the $G$-action on each
$U_j$ is free and the real codimension of $V_j - U_j$ in $V_j$ is more
than $2n-i+1$.  Then $V_1 \oplus V_2$ contains an open set $W$ on
which $G$ acts freely and which contains both $U_1 \oplus V_2$ and
$V_1 \oplus U_2$.  The dimension of
\[
(X \times W)/G - (X \times (U_1 \oplus V_2))/G
\]
is less than $2l_1 + 2l_2 - 2g+i-1$, so
\[
H_{i+2l_1+2l_2-2g,BM}((X \times W)/G) =
H_{i+2l_1+2l_2-2g,BM}((X \times (U_1 \oplus V_2))/G)
\]
by ({\bf BM6}) and ({\bf BM7}).  But $(X \times (U_1 \oplus V_2))/G$ is
a vector bundle of rank $l_2$ over $(X \times U_1)/G$, so
\[
H_{i+2l_1+2l_2-2g,BM}((X \times W)/G) =
H_{i + 2l_1 - 2g,BM}((X \times U_1)/G)
\]
by ({\bf BM9}).  Similarly,
\[
H_{i+2l_1+2l_2-2g,BM}((X \times W)/G) =
H_{i + 2l_2 - 2g,BM}((X \times U_2)/G).
\]

If the real codimension of the open set $U$ in the representation $V$
is $c$ then $\pi_j(U) = 0$ for $j<c-1$, so the mixed spaces $X_G$ are
algebraic approximations to the Borel space $(X \times EG)/G$.
Combining the construction above with the discussion in section A1
immediately yields\footnote{This is entirely parallel to section 2.3
  of \cite{Edidin--Graham:intersection}.} the following properties:
\begin{itemize}
\item[{\bf EBM1}] There are {\em cap products}
  \[
  H^j_G(X) \otimes H^G_{k,BM}(X) \to H^G_{k-j,BM}(X).
  \]
\item[{\bf EBM2}] If $X$ is a smooth variety of dimension $n$ then
  there is a {\em Poincar\'e duality} isomorphism
  \[
  H^k_G(X) \to H_{2n-k,BM}^G(X).
  \]
\item[{\bf EBM3}] If the action of $G$ on $X$ is trivial then
  \[
  H_{k,BM}^G(X) = \bigoplus_{i+j=k} H_{i,BM}(X) \otimes
  H_{j,BM}^G(pt).
  \]
\item[{\bf EBM4}] There are {\em covariant push-forwards} for proper
  $G$-equivariant maps $f:X \to Y$,
  \[
  f_\star:H^G_{k,BM}(X) \to H^G_{k,BM}(Y).
  \]
\item[{\bf EBM5}] There are {\em contravariant pull-backs} for
  $G$-equivariant open embeddings $j:U \to Y$,
  \[
  j^\star:H^G_{k,BM}(Y) \to H^G_{k,BM}(U).
  \]
\item[{\bf EBM6}] There is a {\em long exact sequence}
  \[
  \ldots \to H^G_{i+1,BM}(U) \to H^G_{i,BM}(X) \xrightarrow{i_\star}
  H^G_{i,BM}(Y) \xrightarrow{j^\star} H^G_{i,BM}(U) \to \ldots
  \]
  where $j:U \to Y$ is a $G$-equivariant open embedding and $i:X \to
  Y$ is the $G$-equivariant closed embedding of the complement $X$ to
  $U$ in $Y$.
\item[{\bf EBM7}] We have $H^G_{i,BM}(X)=0$ for $i>2n$.
\item[{\bf EBM8}] There is a {\em cycle map}
  \[
  \cl:A^G_k(X) \to H^G_{2k,BM}(X)
  \]
  which is covariant for proper maps and compatible with
  $G$-equivariant Chern classes.
\item[{\bf EBM9}] There are {\em Gysin maps}
  \[
  f^\star:H^G_{k,BM}(X) \to H^G_{k-2d, BM}(Y)
  \]
  for $G$-equivariant l.c.i. morphisms $f:Y \to X$ of relative
  dimension $d$. These are functorial and compatible with the cycle
  class.  When $Y$ is a $G$-equivariant vector bundle over $X$ of rank
  $d$, $f^\star$ is the {\em Thom isomorphism} $H^G_{k,BM}(X) \to
  H^G_{k+2d,BM}(Y)$.
\end{itemize}

\subsection*{A3. Borel--Moore homology groups for quotient stacks}

In this section, we define ordinary and $\Cstar$-equivariant
Borel--Moore homology groups for certain quotient stacks, following
\cite[section 5]{Edidin--Graham:intersection} and \cite[Appendix
C]{Graber--Pandharipande}.  This allows us to consider the
$\Cstar$-equivariant Borel--Moore homology of moduli spaces of stable maps.

\subsubsection*{Non-equivariant Borel--Moore homology}

Given a quotient stack of the form $[X/G]$, where $X$ is a scheme with
a reasonable action of the $g$-dimensional linear algebraic group $G$,
we define the Borel--Moore homology groups of $[X/G]$ to be
\[
H_{i,BM}([X/G]) := H_{i+2g,BM}^G(X).
\]
We can see that this is well-defined using the argument of
\cite[Proposition 16]{Edidin--Graham:localization}.  Suppose that
$[X/G] \cong [Y/H]$ as quotient stacks, where $G$ (respectively $H$)
acts reasonably on the scheme $X$ (respectively $Y$).  Let $V_1$ be an
$l_1$-dimensional representation of $G$ containing an open set $U_1$
on which the $G$-action is free, and let $X_G = (X \times U_1)/G$.
Let $V_2$ be an $l_2$-dimensional representation of $H$ containing an
open set $U_2$ on which the $H$-action is free, and let $Y_H = (Y
\times U_2)/H$.  The diagonal of a quotient stack is representable, so
the fiber product
\[
Z = X_G \times_{[X/G]} Y_H
\]
is a scheme.  But $Z$ fibers over $X_G$ with fiber $U_2$ and over
$Y_H$ with fiber $U_1$, so
\[
H_{i+2l_1,BM}(X_G) = H_{i+2l_1+2l_2,BM}(Z) = H_{i+2l_2,BM}(Y_H).
\]

\subsubsection*{$\Cstar$-equivariant Borel--Moore homology} Here we
follow Appendix C of \cite{Graber--Pandharipande}.  We define the
$\Cstar$-equivariant Borel--Moore homology groups of a quotient stack
$X$ by setting
\begin{equation}
  \label{eq:equivariant_BM_stacks}
  H^{\Cstar}_{i,BM}(X) := H_{i+2l-2,BM}([(X \times U) / \Cstar])
\end{equation}
where $U$ is an open set in an $l$-dimensional representation of
$\Cstar$ as above.  In other words, we follow the prescription
described in section A2 but construct the mixed space $X_{\Cstar}$ as
a stack quotient.  In the case where $X$ is the quotient of a scheme
$Y$ by a reasonable and proper action of a linear algebraic group $G$
such that the $\Cstar$-action on $X$ descends from a reasonable action
of $G \times \Cstar$ on $Y$, we can use the constructions described
earlier in this section to define the right-hand side of
\eqref{eq:equivariant_BM_stacks}.  In applications to moduli stacks of
stable maps, we need only consider quotients of this form where $G =
PGL$ \cite{Fulton--Pandharipande}.

\subsection*{A4. Localization in $\Cstar$-equivariant Borel--Moore
  homology}

This section contains the proof of the localization theorem which we
used in section \ref{section:push-forward}.  In summary: the argument
given by Graber and Pandharipande in Appendix C of
\cite{Graber--Pandharipande} works for Borel--Moore homology too.

\begin{theorem*}
  Suppose that the stack $X$ is the quotient of a scheme $Y$ by a
  reasonable and proper action of a connected reductive group $G$, and
  that $X$ is equipped with a $\Cstar$-action which descends from a
  reasonable action of $G \times \Cstar$ on $Y$.  Then the
  push-forward
  \[
  i_\star: H_{\bullet,BM}^\Cstar(X^\Cstar) \to H^\Cstar_{\bullet,BM}(X)
  \]
  along the inclusion $i:X^\Cstar \to X$ of the $\Cstar$-fixed stack
  becomes an isomorphism after tensoring with the field of fractions
  $\CC(z)$ of $H^\bullet_{S^1}(pt)$.
\end{theorem*}

\begin{proof}
  In view of ({\bf EBM6}) if suffices to show that the
  $\Cstar$-equivariant Borel--Moore homology groups of $X - X^\Cstar$
  vanish after localization.  But $\Cstar$ acts without fixed points
  on $X - X^\Cstar$, so $X - X^\Cstar$ is the quotient of a scheme $Z$
  by a reasonable and proper action of $G \times \Cstar$ and
  \[
  H^\Cstar_{\bullet,BM}(X - X^\Cstar) =  H_{\bullet,BM}([Z/(G \times
  \Cstar)]).
  \]
  But these groups are non-zero in only finitely many degrees, since
  they are isomorphic to Borel--Moore homology groups of the coarse
  quotient.  They therefore vanish after localization.
\end{proof}

\end{document}